\documentclass[12pt]{article}
\def\to{\rightarrow}
\newtheorem{theorem}{Theorem}
\begin{document}
\title{Asymptotic expansions for ratios  of products of  gamma functions}
\author {Wolfgang B\"uhring \\
Physikalisches Institut, Universit\"at Heidelberg, Philosophenweg~12\\ 69120 Heidelberg,
 Germany}
\bigskip
\date{}
\maketitle

\bigskip
buehring@physi.uni-heidelberg.de

\bigskip
\begin{abstract}
 An asymptotic expansion for a ratio of products of gamma functions is derived.  \end{abstract}

\bigskip
2000 Mathematics Subject Classification:  Primary 33B15; Secondary 33C20

\bigskip
Keywords and phrases: Gamma function,  generalized hypergeometric functions

\section{Introduction}

An asymptotic expansion for a ratio of products of gamma functions has recently been found \cite{b00}, which, with
\begin{equation} 
s_1=b_1-a_1-a_2,
\label{e1}
\end{equation}
may be written
\begin{equation}
{{\Gamma (a_1+n)\Gamma (a_2+n)} \over {\Gamma (b_1+n)\Gamma (-s_1+n)}}=1+\sum\limits_{m=1}^M {{{(s_1+a_1)_m(s_1+a_2)_m} \over {(1)_m(1+s_1-n)_m}}}+O(n^{-M-1})
\label{e2}
\end{equation}
as 
$n\to \infty $. Here use is made of the Pochhammer symbol
\[ (x)_n=x(x+1)\cdots (x+n-1)=\Gamma (x+n)/\Gamma (x). \]

The special case when $b_1=1$ of this formula (2) had been stated earlier by Dingle \cite{din}, and there were proofs by Paris \cite{par} and Olver \cite{olv1}\cite{olv2}.

The proof of (\ref{e2}) is based on the formula for the analytic continuation near unit argument of the Gaussian hypergeometric function. For the more general hypergeometric functions
\begin{equation}
_{p+1}F_p\,\left( \left.\begin{array}{c}a_1,a_2,\ldots ,a_{p+1}\\
  b_1,\ldots ,b_p\end{array} \right| z \right)=\sum\limits_{n=0}^\infty  {{{(a_1)_n (a_2)_n\cdots (a_{p+1})_n} \over {(b_1)_n\cdots (b_p)_n (1)_n }}z^n},\quad (|z|<1),
\label{e3}
\end{equation}
the  analytic continuation near $z=1$ is known too, and this raises the question as to whether a sufficiently simple asymptotic expansion can be derived in a similar way for a ratio of products of more gamma function factors. This is indeed the case, and it is the purpose of this work to present such an expansion.

\section{Derivation of the asymptotic expansion}

The analytic continuation of the hypergeometric function near unit argument may be written \cite{b92}
\begin{equation}
{{\Gamma (a_1)\Gamma (a_2)\cdots \Gamma (a_{p+1})} \over {\Gamma (b_1)\cdots \Gamma (b_p)}}\;_{p+1}F_p\left( \left.\begin{array}{c}a_1,a_2,\ldots ,a_{p+1}\\
  b_1,\ldots,b_p\end{array} \right| z \right)
\label{e4}
\end{equation}
\[=\sum\limits_{m=0}^\infty  {g_m(0)\,(1-z)^m}+(1-z)^{s_p}\sum\limits_{m=0}^\infty  {g_m(s_p)\,(1-z)^m},
\]
where
\begin{equation}
s_p=b_1+\cdots +b_p-a_1-a_2-\cdots -a_{p+1}
\label{e5}
\end{equation}
and the coefficients $g_m$ are known. While the $g_m(0)$ are not needed for the present purpose, the $g_m(s_p)$ are \cite{b92}
\begin{equation}
g_m(s_p)=(-1)^m{{(a_1+s_p)_m(a_2+s_p)_m\Gamma (-s_p-m)} \over {(1)_m}}
\label{e6}
\end{equation}
\[\times \sum\limits_{k=0}^m  {{{(-m)_k} \over {(a_1+s_p)_k(a_2+s_p)_k}}A_k^{(p)}},\]
where the coefficients $A_k^{(p)}$ will be shown below. 

The left-hand side $L$ of (\ref{e4}) is
\begin{equation}
L=\sum\limits_{n=0}^\infty  {{{\Gamma (a_1+n)\Gamma (a_2+n)\cdots \Gamma (a_{p+1}+n)} \over {\Gamma (b_1+n)\cdots \Gamma (b_p+n)\Gamma (1+n)}}z^n}.
\label{e7}
\end{equation}
The asymptotic behaviour, as $n\to \infty $, of the coefficients of this power series is governed \cite{fo}\cite{olv0}\cite{ss} by the terms $R$ on the right-hand side which, at $z=1$, are singular,
\begin{equation}
R=\sum\limits_{m=0}^\infty  {g_m(s_p)(1-z)^{s_p+m}}.
\label{e8}
\end{equation}
Expanded by means of the binomial theorem in its hypergeometric-series-form, this is
\begin{equation}
R=\sum\limits_{m=0}^\infty  {g_m(s_p)\sum\limits_{n=0}^\infty  {{{(-s_p-m)_n} \over {\Gamma(1+n)}}z^n}}.
\label{e9}
\end{equation}
Interchanging the order of summation (and making use of the reflection formula of the gamma function ) we may get
\begin{equation}
R=\sum\limits_{n=0}^\infty  {{{\Gamma (-s_p+n)} \over {\Gamma (1+n)}}\sum\limits_{m=0}^\infty  {(-1)^mg_m(s_p){1 \over {\Gamma (-s_p-m)(1+s_p-n)_m}}}z^n}.
\label{e10}
\end{equation}
Comparison of the coefficients of the two power series for $R$ and $L$, which asymptotically, as $n \to \infty$, should agree, then leads to 
\begin{equation}
{{\Gamma (a_1+n)\Gamma (a_2+n)\cdots \Gamma (a_{p+1}+n)} \over {\Gamma (b_1+n)\cdots \Gamma (b_p+n)\Gamma (-s_p+n)}}
\label{e11}
\end{equation}
\[\sim\sum\limits_{m=0}^\infty  {(-1)^mg_m(s_p){1 \over {\Gamma (-s_p-m)(1+s_p-n)_m}}}.\]
Inserting  $g_m$ from (\ref{e6}) and keeping the first $M+1$ terms of the asymptotic series, we get 
\begin {theorem}
\begin{equation}
{{\Gamma (a_1+n)\Gamma (a_2+n)\cdots \Gamma (a_{p+1}+n)} \over {\Gamma (b_1+n)\cdots \Gamma (b_p+n)\Gamma (-s_p+n)}}=1
\label{e12}
\end{equation}
\[+\sum\limits_{m=1}^M {{{(a_1+s_p)_m(a_2+s_p)_m} \over { (1)_m(1+s_p-n)_m}}\sum\limits_{k=0}^m {{{(-m)_k} \over {(a_1+s_p)_k(a_2+s_p)_k}}A_k^{(p)}}}+O(n^{-M-1})\]
as
$n\to \infty $ , where
$s_p=b_1+\cdots +b_p-a_1-a_2-\cdots -a_{p+1}$ .
\end{theorem}
The simple formula (\ref{e2}) above corresponding to $p=1$ can be recovered from this theorem if we define $A_0^{(1)}=1$, $A_k^{(1)}=0$ for $k>0$, so that then the sum over $k$ is equal to $1$ and disappears. The coefficients for larger $p$ can be found in \cite{b92}, but a few of them are here displayed again for convenience:
\begin{equation}
A_k^{(2)}={{(b_2-a_3)_k(b_1-a_3)_k} \over {k! }},
\label{e13}
\end{equation}
\begin{equation}
A_k^{(3)}=\sum\limits_{k_2=0}^k {{{(b_3+b_2-a_4-a_3+k_2)_{k-k_2}(b_1-a_3)_{k-k_2}(b_3-a_4)_{k_2}(b_2-a_4)_{k_2}} \over {(k-k_2)! k_2!}}},
\label{e14}
\end{equation}

\begin{equation}
A_k^{(4)}=\sum\limits_{k_2=0}^k {{{(b_4+b_3+b_2-a_5-a_4-a_3+k_2)_{k-k_2}(b_1-a_3)_{k-k_2}} \over {(k-k_2)! }}}
\label{e15}
\end{equation}
  \[\times \sum\limits_{k_3=0}^{k_2} {{{(b_4+b_3-a_5-a_4+k_3)_{k_2-k_3}(b_2-a_4)_{k_2-k_3}} \over {(k_2-k_3)! }}}\]
  \[\times {{(b_4-a_5)_{k_3}(b_3-a_5)_{k_3}} \over {k_3! }}.\]
 For $p=3, 4, ...,$ several other representations are possible \cite{b92}, such like
\begin{equation}
A_k^{(3)}={{(b_3+b_2-a_4-a_3)_k(b_1-a_3)_k} \over {k! }}
\label{e16}
\end{equation}
\[\times \,_3F_2\left( \left. \begin{array}{c}b_3-a_4,b_2-a_4,-k\\
  b_3+b_2-a_4-a_3,1+a_3-b_1-k\end{array} \right|1 \right)\]
or
\begin{equation}
A_k^{(3)}={{(b_1+b_3-a_3-a_4)_k(b_2+b_3-a_3-a_4)_k} \over {k! }}
\label{e17}
\end{equation}
\[\times \,_3F_2\left( \left. \begin{array}{c}b_3-a_3,b_3-a_4,-k\\
  b_1+b_3-a_3-a_4,b_2+b_3-a_3-a_4\end{array} \right|1 \right). \]
For $p=2$, the formula (\ref{e12}) may be written simply as
\begin{equation}
{{\Gamma (a_1+n)\Gamma (a_2+n) \Gamma (a_3+n)} \over {\Gamma (b_1+n) \Gamma (b_2+n)\Gamma (-s_2+n)}}=1
\label{e18}
\end{equation}
\[+\sum\limits_{m=1}^M {{{(a_1+s_2)_m(a_2+s_2)_m} \over { (1)_m(1+s_2-n)_m}}{}_3F_2 \left( \left. \begin{array}{c}b_2-a_3,b_1-a_3,-m\\
  a_1+s_2,a_2+s_2\end{array} \right|1 \right)}+O(n^{-M-1}),\]
where $s_2=b_1+b_2-a_1-a_2-a_3$.

\section{Additional comments}

The derivation of the theorem is based on the continuation formula (\ref{e4}) which holds, as it stands, only if $s_p$ is not equal to an integer. Nevertheless, the theorem is valid without such a restriction. This can be verified if the derivation is repeated starting from any of the continuation formulas for the exceptional cases \cite{b92}. Instead of or in addition to the binomial theorem, the expansion
\begin{equation}
(1-z)^m\ln (1-z)=\sum\limits_{n=1}^\infty  {c_nz^n},
\end{equation}
is then needed for integer $m\ge 0$, where
\begin{equation}
c_n=-{1 \over n}(-1)^m{{\Gamma(1+m)\Gamma (n-m)} \over {\Gamma (n)}}
\end{equation}
for $n>m$, while the coefficients are not needed here for $n\le m$.

The theorem has been proved here for any sufficiently large positive integer $n$ only. On the basis of the discussion in \cite{b00}, it can be suspected that the theorem may be theoretically valid (although less useful) in the larger domain of the complex $n$-half-plane $\hbox{Re}(s_p+a_1+a_2-1+n)\ge 0$.

Expansions for ratios of even more general products of gamma functions are treated in a recent monograph by Paris and Kaminski \cite{paka}.

\end{document}